\documentclass[12pt]{amsart}
\usepackage{amsmath}
\usepackage{enumerate}
\usepackage{amssymb}
\usepackage{amsbsy}
\usepackage{amsfonts}
\theoremstyle{plain}
\newtheorem{thm}{Theorem}
\newtheorem{prop}{Proposition}[section]
\newtheorem{lem}[prop]{Lemma}
\newtheorem{thr}[prop]{Theorem}

\theoremstyle{definition}

\newtheorem{ex}[prop]{Example}
\newtheorem{rem}[prop]{Remark}
\newtheorem{dfn}[prop]{Definition}

\newtheorem*{re*}{Remark}
\newtheorem*{ex*}{Example}
\numberwithin{equation}{section}

\newcommand{\sm}{\left(\begin{smallmatrix}}
\newcommand{\esm}{\end{smallmatrix}\right)}

\newcommand{\C}{\mathbb{C}}

\newfont{\FieldFont}{msbm10 scaled\magstep1}

\def\C{\mathbb C}

\begin{document}

\title[The $p$-adic Limit of Weakly Holomorphic Modular Forms]
{$p$-adic Limit of the Fourier Coefficients of Weakly Holomorphic
Modular Forms of Half Integral Weight}

\begin{abstract}
Serre obtained the p-adic limit of the integral Fourier
coefficients of modular forms on $SL_2(\mathbb{Z})$ for
$p=2,3,5,7$. In this paper, we extend the result of Serre to
weakly holomorphic modular forms of half integral weight on
$\Gamma_{0}(4N)$ for $N=1,2,4$. The proof is based on linear
relations among Fourier coefficients of modular forms of half
integral weight. As applications of our main result, we obtain
congruences on   various modular objects, such as those for
Borcherds exponents, for   Fourier coefficients of quotients of
Eisentein series  and for
   Fourier coefficients of  Siegel modular forms on the Maass Space.
\end{abstract}

\author{D. Choi }
\address{School of Liberal Arts and Sciences, Korea Aerospace University, 200-1, Hwajeon-dong, Goyang, Gyeonggi, 412-791, Korea}
\email{choija@postech.ac.kr}
\author{Y. Choie }
\address{Department of Mathematics and Pohang Mathematical
Institute\\
POSTECH\\
Pohang, 790--784, Korea} \email{yjc@postech.ac.kr}
\subjclass[2000]{11F11,11F33} \keywords{modular forms, $p$-adic
limit,   Borcherds exponents, Maass space }
\thanks{This work
 was partially supported by KOSEF R01-2003-00011596-0 , ITRC and BRSI-POSTECH}

 \maketitle

\today

\section{\bf Introduction and Statement of Main Results}

Serre obtained the p-adic limits of the integral Fourier
coefficients of modular forms on $SL_2(\mathbb{Z})$ for
$p=2,3,5,7$ (see Th\'{e}or\`{e}me 7 and Lemma 8 in \cite{S}). In
this paper, we extend the result of Serre to weakly holomorphic
modular forms of half integral weight on $\Gamma_{0}(4N)$ for
$N=1,2,4$. The proof is based on linear relations among Fourier
coefficients of modular forms of half integral weight. As
applications of our main result, we obtain congruences for various
modular objects, such as those for Borcherds exponents, for
Fourier coefficients of quotients of Eisentein series and for
Fourier coefficients of Siegel modular forms on the Maass Space.

For odd $d$, let
$$\left<\left(\begin{smallmatrix}
  1 & h_t \\
  0 & 1
\end{smallmatrix}\right)\right>:=\gamma_t \Gamma_0(4N)_t
\gamma_t^{-1},$$ where $\gamma_t=\left(\begin{smallmatrix}
  a & b \\
  c & d
\end{smallmatrix}\right) \in \Gamma(1)$ and
$\gamma_t(t)=\infty$. We denote the $q$-expansion of a modular form
$f\in M_{\lambda+\frac{1}{2}}(\Gamma_0(4N))$ at each cusp $t$ of
$\Gamma_0(4N)$ by
\begin{equation}\label{cusp}
(f \mid_{\lambda+\frac{1}{2}}
\gamma_t)(z)=(cz+d)^{-\lambda-\frac{1}{2}}f\left(\frac{az+b}{cz+d}\right)
=q_t^{r(t)}\sum_{n=b_t}^{\infty}a_f^t(n)q_t^{n}, \; q_t:=q^{{\frac{1
 }{h_t}}},
\end{equation}
where \begin{equation}\label{1.2} r(t) \in
\left\{0,\frac{1}{4},\frac{1}{2},\frac{3}{4}\right\}.
\end{equation}
When $t \sim \infty$, we denote $a_f^{t}(n)$ by $a_f(n)$. Note
that the number $r(t)$ is independent of the choice of $f\in
M_{\lambda+\frac{1}{2}}(\Gamma_0(4N))$ and $\lambda. $  We call
$t$  a regular cusp if $r(t)=0$ (see  Chapter IV. $\S$1. of
\cite{K} for a more general definition   of a $\lambda$-regular
cusp
  ).

\begin{rem}
Our definition of a regular cusp is different from the usual one.
\end{rem}

Let $U_{4N}:=\{t_{1},\cdots,t_{\nu\left(4N\right)}\}$ be the set
of all inequivalent regular cusps of $\Gamma_0(4N)$. Note that the
genus of $\Gamma_0(4N)$ is zero if and only if $1\leq N\leq 4$.
Let $\mathcal{M}_{\lambda+\frac{1}{2}}(\Gamma_0(4N))$ be the space
of weakly holomorphic modular forms of weight
$\lambda+\frac{1}{2}$ on $\Gamma_0(4N)$ and let
$\mathcal{M}^0_{\lambda+\frac{1}{2}}(\Gamma_0(N))$ denote the set
of $f(z) \in \mathcal{M}_{\lambda+\frac{1}{2}}(\Gamma_0(N))$ such
that the constant term  of its $q$-expansion  at each cusp is
zero. Let $U_p$ be the operator defined  by
$$(f|U_p)(z):=\sum_{n=n_0}^{\infty}a_f(pn)q^n.$$
Let $\mathcal{O}_L$ be the ring of integers of a number field $L$
with a prime ideal $ \mathfrak{p} \subset \mathcal{O}_L$. For
$f(z):=\sum a_f(n)q^n$ and  $g(z):=\sum a_g(n)q^n \in
L[[q^{-1},q]]$ we write $$f(z)\equiv g(z) \pmod{\mathfrak{p}}$$
if and only if $a_f(n)-a_g(n) \in \mathfrak{p}$ for every integer
$n$.

With these notations we state the following theorem.

\begin{thm}\label{main-1} For $N=1,2,4$ consider
$$f(z):=\sum_{n=n_0}^{\infty}a_f(n)q^n \in
\mathcal{M}^0_{\lambda+\frac{1}{2}}(\Gamma_0(4N)) \cap
L[[q^{-1},q]].$$
 Suppose that $ \mathfrak{p} \subset \mathcal{O}_L$
is any prime ideal such that $\mathfrak{p} |p$, $p$ prime, and that
 $a_f(n)$ is $\mathfrak{p}$-integral for every integer $n\geq
 n_0.$

\begin{enumerate}

\item If $p=2$ and $ a_f(0)=0$, then there exists a positive
integer $b$ such that
$$(f |(U_p)^b)(z) \equiv  0 \pmod{\mathfrak{p}^j}  \mbox{\, for each $j\in \mathbb{N}.$}
$$

\item If $p \geq 3$ and $f(z)\in
\mathcal{M}^0_{\lambda+\frac{1}{2}}(\Gamma_0(4N))$ with $\lambda
\equiv 2 \text{ or } 2+\left[\frac{1}{N}\right]
\pmod{\frac{p-1}{2}}$, then there exists a positive integer $b$ such
that
$$(f |(U_p)^b)(z) \equiv 0 \pmod{\mathfrak{p}^j}  \text{ for each } j\in \mathbb{N}.$$
\end{enumerate}
\end{thm}

\begin{rem}
 The $p$-adic limit of a sum of Fourier coefficients of
 $f \in M_{ \frac{3}{2}}(\Gamma_0(4N))$ was studied in \cite{Gu}.
\end{rem}
Our method  only allows to prove a weaker result if $f(z)\not \in
\mathcal{M}^0_{\lambda+\frac{1}{2}}(\Gamma_0(4N))$.
\begin{thm}\label{weak}
For $N=1,2$ or $4$, let
$$f(z):=\sum_{n=n_0}^{\infty}a_f(n)q^n \in
\mathcal{M}_{\lambda+\frac{1}{2}}(\Gamma_0(4N)) \cap
L[[q^{-1},q]].$$ Suppose that $ \mathfrak{p} \subset \mathcal{O}_L$
is any prime ideal with $\mathfrak{p} |p$, $p$ prime, $p\geq 5$, and
that
 $a_f(n)$ is $\mathfrak{p}$-integral for every integer $n\geq n_0$.
If $\lambda \equiv 2 \text{ or } 2+\left[\frac{1}{N}\right]
\pmod{\frac{p-1}{2}}$, then there exists a positive integer $b_0$
such that
\begin{align*}
a_{f}\left(p^{2b-{\bf m}(p:\lambda)}\right)&\equiv -\sum_{t\in
U_{4N}}h_t
a_{\frac{\Delta_{4N,3-\alpha(p:\lambda)}(z)}{R_{4N}(z)^{e\cdot
\omega(4N)}}}^t(0)a^t_f(0) \pmod{\mathfrak{p}}
\end{align*}
for every positive integer $b>b_0$ (see Section 3 for   detailed
notation ).
\end{thm}

\begin{ex}Recall  that the generating function of the overpartition $\bar{P}(n)$
of $n$(see \cite{C-L})
$$\sum_{n=0}^{\infty}\bar{P}(n)q^n=\frac{\eta(2z)}{\eta(z)^2} $$
is in  $\mathcal{M}_{-\frac{1}{2}}(\Gamma_0(16)),$ where
$\eta(z):=q^{\frac{1}{24}}\prod_{n=1}^{\infty}(1-q^n)$. Therefore,
theorem \ref{weak} implies that
$$\bar{P}(5^{2b})\equiv1 \pmod{5}, \forall b \in {\mathbb{N}}.$$

\end{ex}

\section{\bf{Applications: More Congruences  }}

In this section, we study congruences for  various modular objects
such as those for Borcherds exponents and
 for quotients of Eisenstein series.
\subsection{\bf $p$-adic Limits of Borcherds Exponents} Let $\mathcal{M}_H$
denote the set of meromorphic modular forms of integral weight on
$SL_2(\mathbb{Z})$ with Heegner divisor, integer coefficients and
leading coefficient 1. Let
$$\mathcal{M}_{\frac{1}{2}}^+(\Gamma_0(4)):=\{f(z)=\sum_{n=m}^{\infty}a_f(n)q^n
\in \mathcal{M}_{\frac{1}{2}}(\Gamma_0(4))\; | \; a(n)=0 \text{
for } n\equiv 2,3 \pmod{4}\}.$$ If
$f(z)=\sum_{n=n_0}^{\infty}a_f(n)q^n \in
\mathcal{M}_{\frac{1}{2}}^+(\Gamma_0(4))$, then define
$\Psi(f(z))$ by
$$\Psi(f(z)):=q^{-h}\prod_{n=1}^{\infty}(1-q^n)^{a_f(n^2)},$$
where $h=-\frac{1}{12}a_f(0)+\sum_{1<n\equiv 0,1
\pmod{4}}a_f(-n)H(-n).$  Here $H(-n)$ denotes the usual Hurwitz
class number of discriminant $-n$. The following was proved by
Borcherds.

\begin{thr}[\cite{bo}]\label{bor}
The map $\Psi$ is an isomorphism from
$\mathcal{M}_{\frac{1}{2}}^+(\Gamma_0(4))$ to $\mathcal{M}_H$, and
the weight of $\Psi(f(z))$ is $a_f(0)$.
\end{thr}
\noindent Let $j(z)$ be the usual $j$-invariant function with the
product expansion
$$j(z)=q^{-1}\prod_{n=1}^{\infty}(1-q^n)^{A(n)}.$$
Let $F(z):= q^{-h}\prod_{n=1}^{\infty}(1-q^n)^{c(n)}$ be a
meromorphic modular form of weight $k$ in $\mathcal{M}_H$. The
$p$-adic limit of $\sum_{d|n}d\cdot c(d)$ was studied in
\cite{B-O} for $p=2,3,5,7$. Here  we obtain the $p$-adic limit of
$c(d)$ for $p=2,3,5,7$.

\begin{thm}\label{main-3}
Let $F(z):= q^{-h}\prod_{n=1}^{\infty}(1-q^n)^{c(n)}$ be a
meromorphic modular form of weight $k$ in $\mathcal{M}_H$.
\begin{enumerate}
\item If $p=2$, then  for each $j \in \mathbb{N}$  there exists a
positive integer $b$ such that
\[c(mp^b)\equiv 2k \pmod{p^j}\]
for every positive integer $m$. \item If $p\in \{3,5,7\}$, then,
for each $j \in \mathbb{N}$  there exists a positive integer $b$
such that
\[5c(mp^b)-\varpi(F)A(mp^b)\equiv 10k  \pmod{p^j}\]
for every positive integer $m$. Here, $\varpi(F)$ is a constant
determined by the constant term of the $q$-expansion of
$\Psi^{-1}(F)$ at $0$.
\end{enumerate}
\end{thm}

\subsection{Sums of $n$-Squares} For $u \in \mathbb{Z}_{>0},$ let
$$r_{n}(u):=\sharp\{(s_1,\cdots,s_{n}) \in \mathbb{Z}^n \; : \;
s_1^2+\cdots+s_{n}^2=u\}.$$

\begin{thm}\label{main3}
Suppose that $p\geq5$ is a prime. If $\lambda \equiv 2 \text{ or } 3
\pmod{\frac{p-1}{2}}$, then there exists a positive integer $C_0$
such that
$$r_{2\lambda+1}\left(p^{2b-{\bf m}\left(p:\lambda\right)}\right)\equiv  -\left(14-4\alpha\left(p:\lambda\right)\right)
+16\left(\frac{-1}{p}\right)^{\left[\frac{\lambda}{p-1}\right]
+\alpha(p:\lambda){\bf m}(p:\lambda)} \pmod{p},
$$for every $b>C_0.$
\end{thm}
\begin{rem} As for an example, if $\lambda \equiv 2 \pmod{p-1}$ and $p$ is an odd
prime, then there exists a positive integer $C_0$ such that
$$r_{2\lambda+1}\left(p^{2b}\right)\equiv  10 \pmod{p}, \forall b>
C_0$$
\end{rem}

\subsection{Quotients of Eisenstein Series} Congruences for the coefficients of quotients of
elliptic Eisenstein series have been studied in \cite{B-Y}. Let us
consider the Cohen Eisenstein series
$H_{r+\frac{1}{2}}(z):=\sum_{N=0}^{\infty}H(r,N)q^n$  of  weight
$r+\frac{1}{2}, r \geq 2 $ (see \cite{Co}). We derive congruences
for the coefficients of quotients of $H_{r+\frac{1}{2}}(z)$  and
Eisenstein series.

\begin{thm}\label{main5}
Let
$$F(z):=\frac{H_{\frac{5}{2}}(z)}{E_4(z)}=\sum_{n=0}^{\infty}a_{F}(n)q^n,$$
$$G(z):=\frac{H_{\frac{7}{2}}(z)}{E_6(z)}=\sum_{n=0}^{\infty}a_{G}(n)q^n$$
and
$$W(z):=\frac{H_{\frac{9}{2}}(z)}{E_6(z)}=\sum_{n=0}^{\infty}a_{W}(n)q^n.$$
Then there exists a positive integer $C_0$ such that
$$\begin{array}{l}
  a_{F}(11^{2b+1}) \equiv 1 \pmod{11}, \\
  a_{G}(11^{2b+1}) \equiv 6 \pmod{11}, \\
  a_{W}(11^{2b+1}) \equiv 2 \pmod{11}
\end{array}$$
for every integer $b>C_0$.
\end{thm}

\subsection{The Maass Space} Next we deal with congruences for the
Fourier coefficients of a Siegel modular form in the Maass space.
To define the Maass space, let us introduce notations given in
\cite{K-K}: let $T\in M_{2g}(\mathbb{Q})$ be a rational,
half-integral, symmetric, non-degenerate matrix of size $2g$ with
discriminant
$$D_T:=(-1)^g\det(2T).$$
Let  $D_T=D_{T,0}f_T^2,$ where $ D_{T,0}$ is the corresponding
fundamental
discriminant. Furthermore,  let $$G_8:=\left(%
\begin{array}{cccccccc}
  2 & 0 & -1 & 0 & 0 & 0 & 0 & 0 \\
  0 & 2 & 0 & -1 & 0 & 0 & 0 & 0 \\
  -1 & 0 & 2 & -1 & 0 & 0 & 0 & 0 \\
  0 & -1 & -1 & 2 & -1 & 0 & 0 & 0 \\
  0 & 0 & 0 & -1 & 2 & -1 & 0 & 0 \\
  0 & 0 & 0 & 0 & -1 & 2 & -1 & 0 \\
  0 & 0 & 0 & 0 & 0 & -1 & 2 & -1 \\
  0 & 0 & 0 & 0 & 0 & 0 & -1 & 2 \\
\end{array}
\right)$$ and $G_7$ be the upper $(7,7)$-submatrix of $G_8$.
Define
$$S_g:=\left\{\begin{array}{lll}
  G_8^{\bigoplus(g-1)/8}\bigoplus2, &\text{ if } &g\equiv 1 \pmod{8},  \\
   G_8^{\bigoplus(g-7)/8}\bigoplus G_7, &\text{ if } &g\equiv -1 \pmod{8}.  \\
\end{array}\right.$$ For each $m \in \mathbb{N}$  such that $(-1)^gm\equiv 0,1 \pmod{4}$, define
a rational, half-integral, symmetric, positive definite matrix
$T_m$ of size $2g$ by
$$T_m:=\left\{\begin{array}{ll}
  \left(%
\begin{array}{cc}
  \frac{1}{2}S_g & 0 \\
  0 & m/4 \\
\end{array}%
\right), &\text{ if } m\equiv 0 \pmod{4},\\
   \left(%
\begin{array}{cc}
  \frac{1}{2}S_g & \frac{1}{2}e_{2g-1} \\
  \frac{1}{2}e_{2g-1}' & [m+2+(-1)^n]/4 \\
\end{array}%
\right), &\text{ if } m\equiv (-1)^g \pmod{4}  \\
\end{array}
\right.
$$
Here  $e_{2g-1} \in \mathbb{Z}^{(2n-1,1)}$ is the standard column
vector and $e_{2g-1}'$ is its transpose.

\begin{dfn}{\bf (The Maass Space)} Take  $g, k \in \mathbb{N} $ such that
 $g\equiv 0,1 \pmod{4}$ and
 $g\equiv k \pmod{2}. $  Let
\begin{align*}&S_{k+g}^{Maass}(\Gamma_{2g})\\
&\quad:=\left\{F(Z)=\sum_{T>0}A(T)q^{tr(TZ)} \in
S_{k+g}(\Gamma_{2g}) \; \left| \;
A(T)=\sum_{a|f_T}a^{k-1}\phi(a;T)A(T_{|D_T|/a^2})\right\}\right.
\end{align*}
(see (\ref{phi}) for   details). This space is called the Maass
space of genus $2g$ and weight $g+k$.
\end{dfn}
In \cite{K-K} it was proved that the Maass space is the same as
the image of the  Ikeda lifting when $g\equiv 0,1 \pmod{4}$. Using
this fact together with Theorem \ref{main-1}, we derive the
following congruences for  the Fourier coefficients of $F(Z)$ in
$S_{k+g}^{Maass}(\Gamma_{2g})$.

\begin{thm}\label{main6} For  $g\equiv 0,1 \pmod{4},$
let $$F(Z):=\sum_{T>0}A(T)q^{tr(TZ)} \in
S_{k+g}^{Maass}(\Gamma_{2g})$$ with  integral coefficients $A(T)$,
$T>0$. If $k \equiv 2 \text{ or } 3
\pmod{\left.\frac{p-1}{2}\right.}$ for some prime $p$, then, for
each $j\in \mathbb{N}$,  there exists a positive integer $b$ for
which
$$A(T)\equiv 0 \pmod{p^j }$$
for every $T>0, \det(2T)\equiv 0 \pmod{p^{b}}$.
\end{thm}

This paper is organized as follows. Section 3 gives a linear
relation among   Fourier coefficients of modular forms of half
integral weight.  The remaining sections contain detailed proofs
of the main theorems.

\section{\bf{Linear Relation among   Fourier Coefficients of modular
forms of Half Integral Weight}}

Let $V(N;k,n)$ be the subspace of $\C^n$ generated by the first
$n$ coefficients of the $q$-expansion of $f$ at $\infty$ for $f
\in S_k(\Gamma_0(N))$, where $S_k(\Gamma_0(N))$ denotes the space
of cusp forms of weight $k\in \mathbb{Z}$ on $\Gamma_0(N)$. Let
$L(N;k,n)$ be the orthogonal complement of $V(N;k,n)$ in $\C^n$
with the usual inner product of $\C^n$. The vector space
$L(1;k,d(k)+1)$, $d(k)=\dim(S_k(\Gamma(1)))$, was studied by
Siegel to evaluate  the value of the Dedekind zeta function at a
certain point. The vector space $L(1;k,n)$ is explicitly described
in terms of the principal part of negative weight modular forms in
\cite{C-K-O}. These results were extended in \cite{C} to the
groups $\Gamma_0(N)$ of genus zero.  For $1\leq N\leq 4$, let
\begin{align*}
&EV\left(4N,\lambda+{\frac{1}{2}};n\right)\\
&\qquad:=\left\{\left(a^{t_1}_f(0),\cdots,a_f^{t_{\nu(4N)}}(0),a_f(1),\cdots,a_f(n)\right)
\in \mathbb{C}^{n+\nu(4n)} \; \left| \; f \in
M_{\lambda+{\frac{1}{2}}}(\Gamma_0(4N))\right\}\right.,
\end{align*}
 where
$U_{4N}:=\{t_{1},\cdots,t_{\nu\left(4N\right)}\}$ is the set of
all inequivalent regular cusps of $\Gamma_0(4N)$. We define
$EL(4N,\lambda+{\frac{1}{2}};n)$ to be the orthogonal complement
of $EV(4N,\lambda+\frac{1}{2};n)$ in
$\mathbb{C}^{n+\nu\left(4N\right)}$.

Let
$\Delta_{4N,\lambda}:=q^{\delta_{\lambda}(4N)}+O(q^{\delta_{\lambda}(4N)+1})$
be in $M_{\lambda+\frac{1}{2}}(\Gamma_0(4N)$ with the maximum
order at $\infty$, that is, its order at $\infty$ is bigger than
that of any other modular form of the same level and weight.
Furthermore, let
$$R_{4}(z):=\frac{\eta(4z)^{8}}{\eta(2z)^4}, \;
R_{8}(z):=\frac{\eta(8z)^{8}}{\eta(4z)^{4}},$$
$$R_{12}(z):=\frac{\eta(12z)^{12}\eta(2z)^{2}}{\eta(6z)^{6}\eta(4z)^{4}}
\;\text{and } R_{16}(z):=\frac{\eta(16z)^{8}}{\eta(8z)^{4}}.$$ For
$\ell, \; n \in \mathbb{N}$, define $${\bf
m}(\ell:n):=\left\{\begin{array}{c}
0 \text{ if } \left[\frac{2n}{\ell-1}\right]\equiv 0 \pmod{2} \\
1 \text{ if } \left[\frac{2n}{\ell-1}\right]\equiv 1 \pmod{2}
\end{array}\right.$$ and $$\alpha(\ell:n):=n-\frac{\ell-1}{2}\left[\frac{2n}{\ell-1}\right].$$

Let $\omega(4N)$ be the order of zero of $R_{4N}(z)$ at $\infty$.
Note that $R_{4N}(z)\in M_{2}(\Gamma_0(4N))$ has its only zero at
$\infty$. So, using the definition of
$\eta(z)=q^{\frac{1}{24}}\prod_{n=1}^{\infty}(1-q^n)$, we find
that
\begin{equation}\label{omega}
\omega(4)=1,\omega(8)=2,\omega(12)=4,\omega(16)=4.
\end{equation}
For each $g \in M_{r+\frac{1}{2}}(\Gamma_0(4N))$ and $e\in
\mathbb{N}$, let
\begin{equation}\label{1.3}\frac{g(z)}{R_{4N}(z)^{e}}=\sum_{\nu=1}^{e\cdot\omega(4N)}
b(4N,e,g;\nu)q^{-\nu}+O(1) \mbox{ at } \infty. \end{equation} With
these notations we state the following theorem:
\begin{thr}\label{main1}
Suppose that $\lambda \geq 0$ is an integer and $1\leq N\leq 4$.
For each $e\in \mathbb{N}$ such that $e \geq \frac{\lambda}{2}-1,
$ take $r=2e-\lambda+1$. The linear map
$\Phi_{r,e}(4N):M_{r+\frac{1}{2}}(\Gamma_0(4N))\rightarrow
EL(4N,\lambda+\frac{1}{2};e\cdot\omega(4N))$, defined by
\[\begin{array}{l}
  \Phi_{r,e}(4N)(g) \\
  \quad=\left(h_{t_1}
a_{\frac{g(z)}{R_{4N}(z)^{e}}}^{t_1}(0),\cdots,h_{t_{\nu(4N)}}
a_{\frac{g(z)}{R_{4N}(z)^{e}}}^{t_{\nu(4N)}}(0),b(4N,e,g;1),
\cdots , b(4N,e,g;e\cdot\omega(4N)) \right),
\end{array}
\] is an isomorphism.
\end{thr}
\begin{proof}[Proof of Theorem \ref{main1}]
Suppose that $G(z)$ is a meromorphic modular form of weight $2$ on
$\Gamma_0(4N)$. For $\tau \in \mathbb{H} \cup C_{4N}$, let
$D_{\tau}$ be the image of $\tau$ under the canonical map from
$\mathbb{H} \cup C_{4N}$ to a compact Riemann surface $X_0(4N)$.
Here  $\mathbb{H}$ is the usual complex upper half plane, and
$C_{4N}$ denotes the set of all inequivalent cusps of
$\Gamma_0(4N)$. The residue $\textrm{Res}_{D_{\tau}}Gdz$ of $G(z)$
at $D_{\tau} \in X_0(4N)$ is well-defined since we have a
canonical correspondence between a meromorphic modular form of
weight $2$ on $\Gamma_0(4N)$ and a meromorphic 1-form of
$X_0(4N)$. If $Res_{\tau}G$ denotes the residue of $G$ at $\tau$
on $\mathbb{H}$, then
\[\textrm{Res}_{D_{\tau}}Gdz=\frac{1}{l_{\tau}}Res_{\tau}G. \;\;\;\; \]
 Here  $l_{\tau}$ is the
order of the isotropy group at $\tau$. The residue of $G$ at each
cusp $t \in C_{4N}$ is
\begin{equation}\label{Res}
\textrm{Res}_{D_{t}}Gdz=h_t\cdot\frac{a_G^t(0)}{2\pi i}.
\end{equation}%
 Now we give a proof of Theorem \ref{main1}.

To prove Theorem \ref{main1}, take
\[ G(z)=\frac{g(z)}{R_{4N}(z)^e }f(z),\] where $g \in
M_{r+\frac{1}{2}}(\Gamma_0(4N))$ and
$f(z)=\sum_{n=1}^{\infty}a_f(n)q^n \in
M_{\lambda+\frac{1}{2}}(\Gamma_0(4N))$. Note that $G(z)$ is
holomorphic on $\mathbb{H}$. Since $g(z)$, $R_{4N}(z)$ and $f(z)$
are holomorphic and $R_{4N}(z)$ has no zero on $\mathbb{H}$, it is
enough to compute the residues of $G(z)$ only at all inequivalent
cusps to apply the Residue Theorem. The $q$-expansion of
$\frac{g(z)}{R_{4N}(z)^{e}}f(z)$ at $\infty$ is
\begin{align*}
\frac{g(z)}{R_{4N}(z)^e }f(z) &=\left(\sum_{\nu=1}^{e\cdot
\omega(4N)}
b(4N,e,g;\nu)q^{-\nu}+a_{\frac{g(z)}{R_{4N}(z)^{e}}}(0)+O(q)
\right)\left(\sum_{n=0}^{\infty}a_f(n)q^n \right).
\end{align*} Since $R_{4N}(z)$ has no zero at $t\nsim\infty$,  we have
$$\left.\frac{g(z)}{R_{4N}(z)^e }f(z)\right|_2
\gamma_t  =a_{\frac{g(z)}{R_{4N}(z)^{e}}}^{t}(0)a_f(0)+O(q_t).$$
\noindent Further  note that, for an irregular cusp $t$,
$$a_{\frac{g(z)}{R_{4N}(z)^{e}}}^{t}(0)a_f(0)=0.$$ So  the Residue Theorem and
(\ref{Res}) imply that
\begin{equation}\label{dimL}
\sum_{t\in U_{4N}}h_t a_{\frac{g}{R_{4N}^{e\cdot
\omega(4N)}}}^t(0)a^t_f(0)+\sum_{\nu=1}^{e\cdot
\omega(4N)}b(4N,e,g;\nu)a_f(\nu)=0.
\end{equation}
This shows that $\Phi_{r,e}(4N)$ is well-defined. The linearity of
the map $\Phi_{r,e}(4N)$ is clear.

It remains to check that $\Phi_{r,e}(4N)$ is an isomorphism. Since
there exists no holomorphic modular form of negative weight except
the   zero function, we obtain the injectivity of
$\Phi_{r,e}(4N)$. Note that for $e \geq \frac{\lambda-1}{2},$
\begin{align*}
\dim_{\mathbb{C}}\left(EL\left(4N;\lambda+\frac{1}{2},e \cdot
\omega(4N)\right)\right) &=e\cdot \omega(4N)
+\nu(4N)-\dim_{\mathbb{C}}\left(M_{\lambda+\frac{1}{2}}(\Gamma_0(4N))\right).
\end{align*}
However, the set $C_{4N}, 1 \leq N \leq 4,$ of all inequivalent
cusps of $\Gamma_0(4N)$ are
$$
\begin{array}{l}
C_{4}=\left\{\infty,0,\frac{1}{2}\right\},\\
C_{8}=\left\{\infty,0,\frac{1}{2},\frac{1}{4}\right\},\\
C_{12}=\left\{\infty,0,\frac{1}{2},\frac{1}{3},\frac{1}{4},\frac{1}{6}\right\},\\
C_{16}=\left\{\infty,0,\frac{1}{2},\frac{1}{4},\frac{3}{4},\frac{1}{8}\right\}
\end{array}$$
and it can be checked that
\begin{equation}\label{regularc}
\nu(4)=2,\nu(8)=3,\nu(12)=4,\nu(16)=6
\end{equation}
(see $\S$1 of Chapter 4.  in \cite{K} for details). The dimension
formula of $M_{\lambda+\frac{1}{2}}(\Gamma_0(4N))$ (see Table 1)
together with the results
  in (\ref{omega}) and (\ref{regularc}), implies that
\begin{align*}
\dim_{\mathbb{C}}\left(EL\left(4N,\lambda+\frac{1}{2};e\cdot
\omega(N)\right)\right)
=\dim_{\mathbb{C}}(M_{r+\frac{1}{2}}(\Gamma_0(4N)))
\end{align*}
since $r=2e-\lambda+1$.

\begin{table}[htbp]\label{dimf}
\caption{Dimension Formula for $M_k(\Gamma_0(4N))$}
\begin{tabular}[]{|c|c|c|c|}
\hline $\quad \quad N \quad \quad $ & $\quad k=2n+\frac{1}{2}\quad $ & $\quad k=2n+\frac{3}{2}\quad $  & $\quad k=2n\quad $\\
\hline $\quad N=1 \quad$& $\quad n+1 \quad$ & $\quad n+1 \quad$ & $\quad n+1 \quad$\\
\hline $\quad N=2 \quad$ & $\quad 2n+1 \quad$ & $\quad 2n+2 \quad$ & $\quad 2n+1 \quad$\\
\hline $\quad N=3 \quad$ & $\quad 4n+1 \quad$ & $\quad 4n+3 \quad$ & $\quad 4n+1 \quad$\\
\hline $\quad N=4 \quad$ & $\quad 4n+2 \quad$ & $\quad 4n+4 \quad$ & $\quad 4n+1 \quad$\\
\hline
\end{tabular}
\end{table}
So  $\Phi_{r,e}(4N)$ is surjective since the map $\Phi_{r,e}(4N)$
is injective. This completes our claim.
\end{proof}
\section{\bf Proofs of Theorem \ref{main-1} and \ref{weak}}

\subsection{Proof of Theorem \ref{main-1}}
First, we obtain linear relations among   Fourier coefficients of
modular forms of half integral weight modulo $\mathfrak{p}$. Let
$$\mathcal{O}_{\mathfrak{p}}:=\{ \alpha \in L \;|\; \alpha
\text{ is } \mathfrak{p}\text{-integral}\}.$$ Let
\begin{align*}
  \widetilde{M}_{\lambda+\frac{1}{2},\;\mathfrak{p}}(\Gamma_0(4N)):=\{
H(z)&=\sum_{n=0}^{\infty}a_H(n)q^n \in
\mathcal{O}_{\mathfrak{p}}/\mathfrak{p}\mathcal{O}_{\mathfrak{p}}[[q^{-1},q]] \; |\\
  & \;H \equiv h \; \pmod{\mathfrak{p}} \mbox{ for some }
h \in\mathcal{O}_{\mathfrak{p}}[[q^{-1},q]] \cap
M_{\lambda+\frac{1}{2}}(\Gamma_0(4N))\}.
\end{align*}
and
\begin{align*}
  \widetilde{S}_{\lambda+\frac{1}{2},\;\mathfrak{p}}(\Gamma_0(4N)):=\{
H(z)&=\sum_{n=0}^{\infty}a_H(n)q^n \in
\mathcal{O}_{\mathfrak{p}}/\mathfrak{p}\mathcal{O}_{\mathfrak{p}}[[q^{-1},q]] \; |\\
  & \;H \equiv h \; \pmod{\mathfrak{p}} \mbox{ for some }
h \in\mathcal{O}_{\mathfrak{p}}[[q^{-1},q]] \cap
S_{\lambda+\frac{1}{2}}(\Gamma_0(4N))\}.
\end{align*}
The following lemma gives the dimension of
$\widetilde{M}_{\lambda+\frac{1}{2},\; \mathfrak{p}}(\Gamma_0(4N))$.
\begin{lem}\label{local}
Take  $\lambda \in \mathbb{N},  $  $1\leq N\leq 4$ and a prime $p$
  such that $$\left\{\begin{array}{ll}
  p \geq 3 &\text{ if } N=1,2,4, \\
  p \geq 5 &\text{ if } N=3. \\
\end{array} \right.$$
Now  take any prime ideal  $ \mathfrak{p} \subset \mathcal{O}_L, $
$\mathfrak{p} |p$. Then
\begin{equation*}
\dim\widetilde{M}_{\lambda+\frac{1}{2},\;\mathfrak{p}}(\Gamma_0(4N))=\dim
M_{\lambda+\frac{1}{2}}(\Gamma_0(4N))
\end{equation*}
and
\begin{equation*}
\dim\widetilde{S}_{\lambda+\frac{1}{2},\;\mathfrak{p}}(\Gamma_0(4N))=\dim
S_{\lambda+\frac{1}{2}}(\Gamma_0(4N)).
\end{equation*}
\end{lem}
\begin{proof} Let
$$j_{4N}(z)=q^{-1}+O(q)$$
be a meromorphic modular function with a pole only at $\infty$.
Explicitly, these functions  are \begin{align*}
  j_4(z)=&\frac{\eta(z)^{8}}{\eta(4z)^{8}}+8,  &j_8(z)=\frac{\eta(4z)^{12}}{\eta(2z)^{4}\eta(8z)^{8}}, \quad \quad \quad  \\
  j_{12}(z)=&\frac{\eta(4z)^{4}\eta(6z)^{2}}{\eta(2z)^{2}\eta(12z)^{4}},
  &j_{16}(z)=\frac{\eta^{2}(z)\eta(8z)}{\eta(2z)\eta^{2}(16z)}+2.
  \end{align*}
Since the Fourier coefficients of $\eta(z)$ and $\frac{1}{\eta(z)}$
are integral, the $q$-expansion of $j_{4N}(z)$ has integral
coefficients.

Recall that
$\Delta_{4N,\lambda}=q^{\delta_{\lambda}(4N)}+O(q^{\delta_{\lambda}(4N)+1})$
is the modular form of weight $\lambda+\frac{1}{2}$ on
$\Gamma_0(4N)$ such that the order of its zero at $\infty$ is
higher  than that of any other modular form  of the same level and
weight. Denote the order of zero of $\Delta_{4N,\lambda}$ at
$\infty$ by $\delta_{\lambda}(4N)$. Then the basis of
$M_{\lambda+\frac{1}{2}}(\Gamma_0(4N))$ can be chosen as
\begin{equation}\label{gener}\left\{\Delta_{4N,\lambda}(z)j_{4N}(z)^e\;\; \left|\;\; 0
\leq e \leq \delta_{\lambda}(4N)\right\}\right..\end{equation} If
$\Delta_{4N,\lambda}(z)$ is $\mathfrak{p}$-integral, then
$\{\Delta_{4N,\lambda}(z)j_{4N}(z)^e  \; | \; 0 \leq e \leq
\delta_{\lambda}(4N)\}$ also forms a basis of
$\widetilde{M}_{\lambda+\frac{1}{2},\mathfrak{p}}(\Gamma_0(4N))$.
Note that $\delta_{\lambda}(4N)=\dim
M_{\lambda+\frac{1}{2}}(\Gamma_0(4N))-1$. So  from Table 1 we have
\begin{equation}\label{delta}
\Delta_{4N,\lambda}(z)=\Delta_{4N,j}(z)R_{4N}(z)^{\frac{\lambda-j}{2}},
\end{equation}
where $\lambda \equiv j \pmod{2}, j \in \{0,1\}$. More precisely,
one can choose $\Delta_{4N,j}(z)$ as  followings:
$$\begin{array}{l}
  \Delta_{4,0}(z)=\theta(z),\;\Delta_{4,1}(z)=\theta(z)^3, \\
  \Delta_{8,0}(z)=\theta(z),\; \Delta_{8,1}(z)=\frac{1}{4}
  \left(\theta(z)^3-\theta(z)\theta(2z)^2\right),\\
  \Delta_{12,0}(z)=\theta(z),\; \Delta_{12,1}(z)=\frac{1}{6}\left(\sum_{x,y,z \in
\mathbb{Z}}q^{3x^2+2\left(y^2+z^2+yz\right)}-\sum_{x,y,z \in
\mathbb{Z}}q^{3x^2+4y^2+4z^2+4yz}\right), \\
  \Delta_{16,0}(z)=\frac{1}{2}\left(\theta(z)-\theta(4z)\right),
   \;\Delta_{16,1}(z)=\frac{1}{8}\left(\theta(z)^3-3\theta(z)^2\theta(4z)
+3\theta(z)\theta(4z)^2-\theta(4z)^3\right).
\end{array}
$$
Since $\theta(z)=1+2\sum_{n=1}^{\infty}q^n$, the coefficients of
the $q$-expansion of $\Delta_{4N,j}(z)$, $j \in \{0,1\}$, are
$\mathfrak{p}$-integral. This completes the proof.
\end{proof}

\begin{rem}
The proof of Lemma \ref{local} implies that the spaces of
$M_{\lambda+\frac{1}{2}}(\Gamma_0(4N))$ for $N=1,2,4$ are
generated by eta-quotients since
$\theta(z)=\frac{\eta(2z)^5}{\eta(z)^2\eta(4z)^2}$.
\end{rem}
For $1 \leq N \leq 4$ set
$$\widetilde{V_S}\left(4N,\lambda+{\frac{1}{2}};n\right)
:=\left\{\left(a_f(1),\cdots,a_f(n)\right) \in
\mathbb{F}_{\mathfrak{p}}^{n} \; | \; f \in
\widetilde{S}_{\lambda+{\frac{1}{2}}}(\Gamma_0(4N))\right\},
\mathbb{F}_{\mathfrak{p}}:=\mathcal{O}_{\mathfrak{p}}/\mathfrak{p}\mathcal{O}_{\mathfrak{p}}.$$
  We define
$\widetilde{L_S}(4N,\lambda+\frac{1}{2};n)$ to be the orthogonal
complement of $\widetilde{V_S}(4N,\lambda+\frac{1}{2};n)$ in
$\mathbb{F}_{\mathfrak{p}}^{n}$. Using Lemma \ref{local}, we
obtain the following proposition.
\begin{prop}\label{local-relation}
Suppose that $\lambda$ is a positive integer and $1\leq N\leq 4$.
For each $e\in \mathbb{N}$, $e \geq \frac{\lambda}{2}-1$, take
$r=2e-\lambda+1$. The linear map
$\widetilde{\psi_{r,e}}(4N):\widetilde{M}_{r+\frac{1}{2},\mathfrak{p}}(\Gamma_0(4N))\rightarrow
\widetilde{L_S}(4N,\lambda+\frac{1}{2};e\cdot\omega(4N))$, defined
by
\[\widetilde{\psi_{r,e}}(4N)(g)=\left(b(4N,e,g;1), \cdots ,
b(N,e,g;e\cdot\omega(4N)) \right),
\] is an isomorphism. Here
$b(4N,e,g;\nu)$ is defined in (\ref{1.3}).
\end{prop}
\begin{proof}
Note that $\dim S_{\frac{3}{2}}(4N)=0$ and that $$\dim
S_{\lambda+\frac{1}{2}}(4N)+N+1+\left[\frac{N}{4}\right]=\dim
M_{\lambda+\frac{1}{2}}(4N)$$(see \cite{C-O}). So, from Lemma
\ref{local} and Table 1, it is enough to show that
$\psi_{r,e}(4N)$ is injective. If $g$ is in the kernel of
$\psi_{r,e}(4N)$, then $\frac{g(z)}{R_{4N}(z)^e}\cdot R_{4N}(z)^e
\equiv 0 \pmod{{\mathfrak{p}}}$ by Sturm's formula (see
\cite{St}). So  we have $g(z)\equiv 0\pmod{{\mathfrak{p}}}$ since
$R_{4N}(z)^e \not \equiv  0 \pmod{{\mathfrak{p}}}$. This completes
the proof.
\end{proof}

\begin{thr}\label{main2-1}
Take a prime $p, N=1,2,4$ and
$$f(z):=\sum_{n=n_0}^{\infty}a_f(n)q^n \in
S_{\lambda+\frac{1}{2}}(\Gamma_0(4N)) \cap L[[q]].$$ Suppose that $
\mathfrak{p} \subset \mathcal{O}_L$ is any prime ideal with
$\mathfrak{p} |p$ and that
 $a_f(n)$ is $\mathfrak{p}$-integral for every integer $n\geq n_0$.
If $\lambda \equiv 2 \text{ or } 2+\left[\frac{1}{N}\right]
\pmod{\frac{p-1}{2}}$ or $p=2$, then there exists a positive
integer $b$ such that
$$a_{f}\left(np^{b}\right) \equiv 0 \pmod{{\mathfrak{p}}}, \forall n \in \mathbb{N}.$$
\end{thr}
\begin{proof}[Proof of Theorem \ref{main2-1}]
{\bf i)} First, suppose that $p\geq3$: Take   positive integers
$\ell$ and $b$ such that
\begin{equation}\label{weight}
\frac{3-2\alpha(p:\lambda)}{2}p^{2b}+\left(\lambda+\frac{1}{2}\right){p^{{\bf
m}(p:\lambda)}}+\ell(p-1)=2.
\end{equation}
 Note that if $b$ is large enough, that is,
$b>\log_{p}\left(\frac{2}{3-2\alpha(p:\lambda)}
\left(\lambda+\frac{1}{2}\right){p^{{\bf
m}(p:\lambda)}}-2\right)$, then there exists a positive integer
$\ell$ satisfying (\ref{weight}). Also  note that $a_f^t(0)=0$ for
every cusp $t$ of $\Gamma_0(4N)$ since $f(z)$ is a cusp form. So,
if $r=2e-\alpha(p:\lambda)+1$, then Theorem \ref{main1} implies
that, for   $g(z)\in \widetilde{M}_{r+\frac{1}{2}}(\Gamma_0(4N))$,
$$\sum_{\nu=1}^{e\cdot \omega(4N)}b(4N,e,g;\nu)a_f(\nu p^{2b- {\bf  m} (p:\lambda)})\equiv 0
\pmod{{\mathfrak{p}}},$$ since
\begin{align*}
\left(\frac{g(z)}{R_{4N}(z)^e }\right)^{p^{2b}}&f(z)^{p^{{\bf  m}
(p:\lambda)}}E^{\ell}_{p-1}(z)\\
&\equiv\left(\sum_{\nu=1}^{e\cdot \omega(4N)} b(4N,e,g;\nu)q^{-\nu
p^{2b}}+a_{\frac{g(z)}{R_{4N}(z)^{e}}}(0)+\sum_{n=1}^{\infty}a_{\frac{g(z)}{R_{4N}(z)^{e}}}(n)q^{n
p^{2b}} \right)\\
&\quad\cdot\left(\sum_{n=0}^{\infty}a_f(n)q^{np^{{\bf m}
(p:\lambda)}} \right) \pmod{p}.
\end{align*}
So   Proposition \ref{local-relation} implies that
$$\begin{array}{l}
  \left(a\left(p^{2b- {\bf  m} (p:\lambda)}\right),a\left(2p^{2b-
{\bf  m} (p:\lambda)}\right), \cdots,a\left(e\cdot\omega(4N)p^{2b-
{\bf  m} (p:\lambda)}\right)\right) \\
 \in
\widetilde{V_S}\left(4N,\alpha(p:\lambda)+{\frac{1}{2}};n\right).
\end{array} $$
If $\alpha(p:\lambda)=2$ or $2+\left[\frac{1}{N}\right]$, then
$$\dim
S_{\alpha(p:\lambda)+\frac{1}{2}}(\Gamma_0(4N)) =\dim
\widetilde{V_S}\left(4N,\alpha(p:\lambda)+{\frac{1}{2}};n\right)=0.$$

{\bf ii)}   $p=2$: Note that
$\frac{\Delta_{4N,1}(z)}{R_{4N}(z)}=q^{-1}+O(1)$ for $N=1,2,4$.
So, there exists a polynomial $F(X)\in \mathbb{Z}[X]$ such that
$$F(j_{4N}(z))\frac{\Delta_{4N,1}(z)}{R_{4N}(z)}=q^{-n}+O(1).$$
For an integer $b$, $2^{2^b}>\lambda+2$, let
$$G(z):=\left(F(j_{4N}(z))\frac{\Delta_{4N,1}(z)}{R_{4N}(z)}\right)
^{2^b}f(z)\theta(z)^{2^{1+2b}-2\lambda+3}.$$ Since
$\theta(z)\equiv 1 \pmod{2}$, Theorem \ref{main1} implies that
$a_{f}(2^b\cdot n)\equiv 0 \pmod{{\mathfrak{p}}}$.
\end{proof}

To apply Theorem \ref{main2-1}, we need the following two
propositions.
\begin{prop}[Proposition 3.2 in \cite{T}]\label{tre}
Suppose that $p$ is an odd prime, $k$ and $N$ are integers with
$(N,p)=1$. Let
$$f(z)=\sum a(n)q^n \in
\mathcal{M}_{\lambda+\frac{1}{2}}(\Gamma_0(4N)).$$ Suppose  that
$\xi:=\left(\begin{smallmatrix}
  a & b \\
  cp^2 & d
\end{smallmatrix}\right)$, with $ac>0$. Then there exist  $n_0, h_0 \in
\mathbb{N}$ with  $h_0|N, $   a sequence $\{a_0(n)\}_{n \geq n_0}$
and $r_0 \in \{0,1,2,3\}$ such that
$$(f|U_{p^m}|_{\lambda+\frac{1}{2}}\xi)(z)=\sum_{\begin{smallmatrix}
  n \geq n_0 \\
  4n+r_0\equiv0 \pmod{p^m}
\end{smallmatrix}}a_0(n)q^{\frac{4n+r_0}{4h_0p^m}}, \; \forall m\geq 1 .
$$
\end{prop}

\begin{prop}[Proposition 5.1 in \cite{A-B}]\label{reduction}
Suppose that $p$ is an odd prime such that $p \nmid N$ and
consider
$$g(z)=\sum_{n=1}^{\infty}a(n)q^n \in
S_{\lambda+\frac{1}{2}}(\Gamma_0(4Np^j))\cap L[[q]], \mbox{\, for
each $j\in \mathbb{N}$} .$$  Suppose further that $ \mathfrak{p}
\subset \mathcal{O}_L$ is any prime ideal with $\mathfrak{p} |p$
and that
 $a(n)$ is $\mathfrak{p}$-integral for every integer $n\geq 1$.
Then there exists $G(z) \in
 S_{\lambda'+\frac{1}{2}}(\Gamma_0(4N))\cap\mathcal{O}_L[[q]]$
such that
$$G(z)\equiv g(z) \pmod{{\mathfrak{p}}},$$ where
$\lambda'+\frac{1}{2}=(\lambda+\frac{1}{2})p^{j}+p^e(p-1)$ with $e
\mathbf{N}$ large.

\end{prop}
\begin{rem}Proposition \ref{reduction} was proved for
$p\geq5$ in \cite{A-B}.  One can   check that this holds also for
$p=3$.
\end{rem}
\noindent Now we prove Theorem \ref{main-1}.
\begin{proof}[Proof of Theorem \ref{main-1}]
Take  $$G_p(z):=\left\{\begin{array}{lll}
  \frac{\eta(8z)^{48}}{\eta(16z)^{24}}\in M_{12}(\Gamma_0(16))& \text{ if } p=2,\\
  \frac{\eta(z)^{27}}{\eta(9z)^3}\in M_{12}(\Gamma_0(9))& \text{ if } p=3,\\
  \frac{\eta(4z)^{p^2}}{\eta(4p^2z)}\in M_{\frac{p^2-1}{2}}(\Gamma_0(p^2))& \text{ if } p \geq 5.
\end{array}\right.$$
Using properties of eta-quotients (see \cite{G-H}), note that
$G_p(z)$ vanishes at every cusp of $\Gamma_0(16)$ except $\infty$
if $p=2$, and vanishes at every cusp $\frac{a}{c}$ of
$\Gamma_0(4Np^2)$ with $p^2 \nmid N$ if $p \geq 3$. Thus,
Proposition \ref{tre} implies that there exist positive integers
$\ell,m,k$ such that
$$\left\{\begin{array}{ll}
  (f |U_{p^m})(z) G_p(z)^{\ell}\in
S_{k+\frac{1}{2}}(\Gamma_0(16)) &\text{ if } p=2, \\
  (f |U_{p^m})(z) G_p(z)^{\ell}\in
S_{k+\frac{1}{2}}(\Gamma_0(4p^2 N)) &\text{ if } p\geq3.
\end{array}\right.$$ Note that $k\equiv \lambda \pmod{p-1}$. Using Proposition
\ref{reduction}, we can find $$F(z)\in
S_{k'+\frac{1}{2}}(\Gamma_0(4N))$$ such that
$F(z)\equiv(f(z)|U_{p^m})G_p(z)^{\ell}\equiv (f |U_{p^m})(z)
\pmod{{\mathfrak{p}}}$ and $k'\equiv k\pmod{p-1}$. Theorem
\ref{main2-1} implies that there exists a positive integer $b$
such that $(F|U_{p^{2b}})(z)\equiv 0 \pmod{{\mathfrak{p}}}$. Thus,
we have shown so far   that if $\rho\in {\mathfrak{p}}\setminus
\mathfrak{p}^2$,   all the Fourier coefficients of
$\frac{1}{\rho}\cdot F(z)|U_{p^{m+2b}}$ are
$\mathfrak{p}$-integral. Repeat this argument to complete our
claim.
\end{proof}

\subsection{Proof of Theorem \ref{weak}} Theorem \ref{weak} can be derived
from Theorem \ref{main1} by taking  a special modular form.

\begin{proof}[Proof of Theorem \ref{weak}] Take a positive
integer $\ell$ and a positive even integer $u$ such that
$$\frac{3-2\alpha(p:\lambda)}{2}p^u+\left(\lambda+\frac{1}{2}\right)
{p^{{\bf m}(p:\lambda)}}+\ell(p-1)=2.$$ Let $F(z):=
\left(\frac{\Delta_{4N,3-\alpha(p:\lambda)}(z)}{R_{4N}(z)}\right)^{p^u}$
and $G(z):=E_{p-1}(z)^{\ell}f(z)^{p^{{\bf m}(p:\lambda)}}$. Since
$E_{p-1}(z) \equiv 1 \pmod{p}$, we have
\[F(z)G(z)\equiv\left(
\sum_{n=-1}^{\infty}a_{\frac{\Delta_{4N,3-\alpha(p:\lambda)}(z)}{R_{4N}(z)}}(n)q^{np^u}\right)
\left(\sum_{n=m_{\infty}}^{\infty}a_f(n)q^{n {\bf m}(p:\lambda)}
\right) \pmod{{\mathfrak{p}}}.
\] If   Fourier coefficients of $f(z)$ at each cusp are
$\mathfrak{p}$-integral, then
\begin{align*}
\left((F \cdot G )|_{2} \gamma_t \right) (z)&\equiv
\left(q_t^r\sum_{n=m_t}^{\infty}a^t_F(n)q_t^n\right)\left(q_t^r\sum_{n=0}^{\infty}a^t_G(n)q_t^n\right)\\
&\equiv \left(q_t^r\sum_{n=m_t}^{\infty}a^t_f(n)q_t^n\right)
\left(q_t^{p^u}\sum_{n=0}^{\infty}a_{\frac{\Delta_{4N,3-\alpha(p:\lambda)}(z)}{R_{4N}(z)}}^t(n)q_t^{p^u}\right)
\pmod{{\mathfrak{p}}}
\end{align*}
for $t\nsim \infty$. Since
$$\begin{array}{lll}
   a_{F(z)G(z)}(0)\equiv a_{\frac{\Delta_{4N,3-\alpha(p:\lambda)}(z)}{R_{4N}(z)}}(0)a_f(0)+a_f(p^{u- {\bf m}(p:\lambda) }     )
   &\pmod{{\mathfrak{p}}} &,\\
   a^{t}_{F(z)G(z)}(0)\equiv a_{\frac{\Delta_{4N,3-\alpha(p:\lambda)}(z)}{R_{4N}(z)}}^t(0)a_f^t(0)
    &\pmod{{\mathfrak{p}}}& \text{ for } t\nsim \infty,
  \end{array}
$$ for large $u$, the Residue Theorem implies Theorem \ref{weak} by letting
$u=2b.$  Therefore  it is enough to check a
$\mathfrak{p}$-integral property of   Fourier coefficients of
$f(z)$ at each cusp:
  take  a positive integer $e$ such that
$\Delta(z)^e f(z)$ is a holomorphic modular form, where
$\Delta(z):=q\prod_{n=1}^{\infty}(1-q^{n})^{24}.$     Note that
the $q$-expansions of $j_{4N}(z)$ and $\Delta_{4N,12e+\lambda}(z)$
at each cusp are $p$-integral. Thus (\ref{gener}) implies that
$$\Delta(z)^ef(z)=\sum_{n=0}^{\delta_{12e+\lambda}(4N)}c_nj_{4N}(z)^n\Delta_{4N,12e+\lambda}(z).$$
Moreover, $c_n$ is ${\mathfrak{p}}$-integral since
$$j_{4N}(z)^n\Delta_{4N,12e+\lambda}(z)=
q^{\delta_{12e+\lambda}(4N)-n}+O\left(q^{\delta_{12e+\lambda}(4N)-n+1}\right)$$
and $f(z) \in \mathcal{O}_L[[q,q^{-1}]]$. Note that $p \nmid 4N$
since $1\leq N \leq 4$ and $p\geq 5$ is a prime. So   Fourier
coefficients of $j_{4N}(z)$, $\Delta_{N,12e+\lambda}(z)$ and
$\frac{1}{\Delta(z)}$ at each cusp  are ${\mathfrak{p}}$-integral.
This completes our claim.
\end{proof}

\section{\bf Proof of Theorem  \ref{main-3}}
  Theorem
\ref{main-3} follows from Theorem \ref{main-1} and Theorem
\ref{bor}.
\begin{proof}[Proof of Theorem \ref{main-3}]
Note that $j(z)\in \mathcal{M}_H$. Let
$$g(z):=\Psi^{-1}(j(z)) \text{ and } f(z):=\Psi^{-1}(F(z))=\sum_{n=n_0}^{\infty}a_f(n)q^n.$$
It is known (see $\S$14 in \cite{bo}) that
$$\frac{1}{3}g(z)=\frac{\frac{d}{dz}(\theta(z))E_{10}(4z)}{4 \pi i \Delta(4z)}-\frac{\theta(z)
\frac{d}{dz}(E_{10}(4z))}{80 \pi i
\Delta(4z)}-\frac{152}{5}\theta(z).$$ Since the constant terms of
the $q$-expansions at $\infty$ of $f(z)$, $\theta(z)$ and $g(z)$
are $0,$  $a_{\theta(z)}^0(0)=\frac{1-i}{2}$ and
$a_g^{0}(0)=\frac{1-i}{2}\cdot\frac{456}{5},$ respectively, we
have
$$f(z)-k\theta(z)-\frac{a_{f}^{0}(0)+k(1-i)/2}{a_g^{0}(0)}g(z)\in \mathcal{M}^0_{\frac{1}{2}}(\Gamma_0(4)).$$
Applying Theorem \ref{main-1}, one obtains the result.
\end{proof}

\section{\bf Proofs of Theorem \ref{main3} and \ref{main5}}
We begin with the following proposition.
\begin{prop} \label{main2}
Let $p$ be an odd prime and $$f(z):=\sum_{n=0}^{\infty}a_f(n)q^n
\in M_{\lambda+\frac{1}{2}}(\Gamma_0(4)) \cap
\mathbb{Z}_{p}[[q]].$$ If $\lambda \equiv 2 \text{ or } 3
\pmod{\frac{p-1}{2}}$, then
\begin{align*} &a_{f}\left(p^{2b-{\bf m}(p:\lambda)}\right)\\
&\qquad\equiv -(14-4\alpha(p:\lambda))
a_{f}(0)+2^{8}\left(2^{-1}-2^{-1}i\right)^{p^b(7-2\alpha(p:\lambda))}a_f^0(0)
\pmod{p}
\end{align*}
for every integer $b>\log_{p}\left(\frac{2}{2\alpha(p:\lambda)-3}
\left(\lambda+\frac{1}{2}\right){p^{{\bf
m}(p:\lambda)}}+2\right)$.
\end{prop}

\begin{proof}[Proof of Proposition \ref{main2}]
For $\nu \in \mathbb{Z}_{\geq0}$,
$$\left(\lambda+\frac{1}{2}\right){p^{ {\bf  m} (p:\lambda)}}:=\nu\cdot(p-1)+\alpha(p:\lambda)+\frac{1}{2}.$$
For an integer $b$ with
$$b>\frac{1}{2}\log_p\left(\frac{2}{3-2\alpha(p:\lambda)}\left(\left(\lambda+\frac{1}{2}\right)p^{{\bf
m}(p:\lambda)}-2\right)\right),$$ there exists  an $\ell \in
\mathbb{N}$ such that
$$\frac{3-2\alpha(p:\lambda)}{2}p^{2b}+\left(\lambda+\frac{1}{2}\right)p^{{\bf m}(p:\lambda)}+\ell(p-1)=2,$$
since
$$
 \frac{3-2\alpha(p:\lambda)}{2}p^{2b}+\left(\lambda+\frac{1}{2}\right)p^{{\bf m}(p:\lambda)}-2
 = \frac{3-2\alpha(p:\lambda)}{2}(p^{2b}-1)+\nu(p-1). $$
We have
$$\begin{array}{l}
  F(z)\equiv \sum_{n=0}^{\infty}a_f(n)q^{n{p^{{\bf m}(p:\lambda)}}} \pmod{p},\\
  G(z)\equiv q^{-p^b}+14-4\alpha(p:\lambda)+a_G(1)q+\cdots\pmod{p}.
\end{array}$$ Note that $a_G(n)$ is $p$-integral for every integer $n$.
Moreover, we obtain
$$F(z) G(z) |_2\left(
                  \begin{smallmatrix}
                    0 & -1 \\
                    1 & 0 \\
                  \end{smallmatrix}
                \right)
                \equiv \left(a_f^0(0)+\cdots\right)
                \left(-2^{6p^b}\left(\frac{1}{2}-\frac{i}{2}\right)^{p^b\left(7-2\alpha(p:\lambda)\right)}
                +\cdots\right) \pmod{p},
$$
where $a_f^0(0)$ is given in (\ref{cusp}). Note that
$\left\{\infty,0,\frac{1}{2}\right\}$ is the set of cusps of
$\Gamma_0(4)$, so  Theorem \ref{weak} implies that
$$a_{f}(p^{2b-{\bf m}(p:n)})+(14-4\alpha(p:\lambda))a_f(0)-2^8a_f^0(0)
\left(\frac{1}{2}-\frac{i}{2}\right)^{p^b\left(7-2\alpha(p:\lambda)\right)}
\equiv 0 \pmod{p}.$$ This proves Proposition \ref{main2}.
\end{proof}

\subsection{Proof of Theorem \ref{main3}}Now we prove Theorem
\ref{main3}.
\begin{proof}[Proof of Theorem \ref{main3}]
Take
$$f(z):=\theta^{2\lambda+1}(z)=1+\sum_{\ell=1}^{\infty}r_{2\lambda+1}(\ell)q^{\ell}=\sum_{n=0}^{\infty}a_f(n)q^n.$$
Note that $f(z)\in M_{\lambda+\frac{1}{2}}(\Gamma_0(4))$. Since
$(\theta |_{\frac{1}{2}}\left(\begin{smallmatrix}
  0 & -1 \\
  1 & 0
\end{smallmatrix}\right)
)(z)=\frac{1-i}{2}+O\left(q^{\frac{1}{4}}\right)$, we obtain
$$a_f(0)=1 \text{ and }a_f^0(0)=\left(\frac{1-i}{2}\right)^{2\lambda+1}.$$
Since $\lambda \equiv 2,3 \pmod{\frac{p-1}{2}}$ and
$\left(\frac{1-i}{2} \right)^8=\frac{1}{16}$, we have
\begin{align*}
&\left(\frac{1}{2}-\frac{i}{2}\right)^{p^{2u}(7-2\alpha(p:\lambda))}a_f^0(0)^{p^{{\bf
m}(p:\lambda)}}\\
&\equiv\left(\frac{1}{2}-\frac{i}{2}\right)^{p^{2u}(7-2\alpha(p:\lambda))}
\left(\frac{1}{2}-\frac{i}{2}\right)^{p^{{\bf
m}(p:\lambda)}\left(2\alpha(p:\lambda)+(p-1)\left(2\left[\frac{\lambda}{p-1}\right]+{\bf m}(p:\lambda)\right)+1\right)}\\
&\equiv\left(\frac{1}{2}-\frac{i}{2}\right)^{(7-2\alpha(p:\lambda))(p^{2u}-1)}
\left(\frac{1}{2}-\frac{i}{2}\right)^{8+2(p-1)
\left[\frac{\lambda}{p-1}\right]+{\bf m}(p:\lambda)p^{{\bf m}(p:\lambda)}(p-1)+(p^{{\bf m}(p:\lambda)}-1)(1+2\alpha(p:\lambda))}\\
&\equiv
\left(\frac{1}{2}-\frac{i}{2}\right)^{8+2\left[\frac{\lambda}{p-1}\right](p-1)+2\alpha(p:\lambda)(p^{{\bf
m}(p:\lambda)}-1)} \equiv
\frac{1}{16}\left(\frac{-1}{p}\right)^{\left[\frac{\lambda}{p-1}\right]+\alpha(p:\lambda){\bf
m}(p:\lambda)}\pmod{p},
\end{align*}
for some  $u \in \mathbf{N}.$
 Applying Proposition \ref{main2}, we
obtain the result.
\end{proof}

\subsection{Proof of Theorem \ref{main5}}
Consider the Cohen Eisenstein series
$H_{r+\frac{1}{2}}(z):=\sum_{N=0}^{\infty}H(r,N)q^n$  of weight
$r+\frac{1}{2}$, where $r\geq2$ is an integer. If $(-1)^rN\equiv
0,1 \pmod{4}$, then $H(r,N)=0$. If $N=0$, then
$H(r,0)=\frac{-B_{2r}}{2r}$. If $N$ is a positive integer and
$Df^2=(-1)^rN$, where $D$ is a fundamental discriminant, then
\begin{equation}\label{3-5}
H(r,N)=L(1-r,\chi_D)\sum_{d|f}\mu(d)\chi_D(d)d^{r-1}\sigma_{2r-1}(f/d).
\end{equation}
Here  $\mu(d)$ is  the $M\ddot{o}bius$ function. The following
theorem implies that the Fourier coefficients of
$H_{r+\frac{1}{2}}(z)$ are $p$-integral if $\frac{p-1}{2}\nmid r$.
\begin{thr}[\cite{Cal}]\label{Cal}
Let $D$ be a fundamental discriminant. If $D$ is divisible by at
least two different primes, then $L(1-n,\chi_D)$ is an integer for
every positive integer $n$. If $D=p$, $p>2$, then $L(1-n,\chi_D)$
is an integer for every positive integer $n$ unless
$\gcd(p,1-\chi_D(g)g^n) \neq 1,$ where $g$ is a primitive root
$\pmod{p}$.
\end{thr}
\begin{proof}[Proof of Theorem \ref{main5}]
Note that $E_{10}(z)=E_{4}(z)E_{6}(z)$. So, $E_{10}(z)F(z)$,
$E_{10}(z)G(z)$ and $E_{10}(z)W(z)$ are modular forms of weights,
$8\cdot\frac{1}{2}$, $7\cdot\frac{1}{2}$ and $8\cdot\frac{1}{2}$
respectively. Moreover, the Fourier coefficients of those modular
forms are $11$-integral, since the Fourier coefficients of
$H_{\frac{5}{2}}(z)$, $H_{\frac{7}{2}}(z)$ and
$H_{\frac{9}{2}}(z)$ are 11-integral by Theorem \ref{Cal}. We have
$$
\begin{array}{l}
  E_{10}(z)F(z)=\frac{B_4}{4}+O(q),\\
  E_{10}(z)F(z)|_{\frac{17}{2}}\left(\begin{smallmatrix}
  0 & -1 \\
  1 & 0
\end{smallmatrix}\right)=\frac{B_4}{4}(1+i)(2i)^{-5}+O\left(q^{\frac{1}{4}}\right), \\
  E_{10}(z)G(z)=\frac{B_6}{6}+O(q), \\
  E_{10}(z)G(z)|_{\frac{15}{2}}\left(\begin{smallmatrix}
  0 & -1 \\
  1 & 0
\end{smallmatrix}\right)=\frac{B_6}{6}(1-i)(2i)^{-7}+O\left(q^{\frac{1}{4}}\right),  \\
  E_{10}(z)W(z)=\frac{B_8}{8}+O(q), \\
  E_{10}(z)W(z)|_{\frac{17}{2}}\left(\begin{smallmatrix}
  0 & -1 \\
  1 & 0
\end{smallmatrix}\right)=\frac{B_8}{8}(1+i)(2i)^{-9}+O\left(q^{\frac{1}{4}}\right),
\end{array}
$$
where $B_{2r}$ is the $2r$th Bernoulli number. The conclusion now
follows from  Proposition \ref{main2}.
\end{proof}
\subsection{\bf Proof of Theorem \ref{main6}} We begin by
introducing some notations (see \cite{K-K}). Let
$V:=(\mathbb{F}_{p}^{2n},Q)$ be the quadratic space over
$\mathbb{F}_p$, where $Q$ is the  quadratic form obtained from a
quadratic form $x\mapsto T[x](x\in \mathbb{Z}_{p}^{2n})$ by
reducing modulo $p$. We denote by $<x,y>:=Q(x,y)-Q(x)-Q(y), \; x,y
\in\mathbb{F}_{p}^{2n}$, the associated bilinear form and let
$$R(V):=\{x\in\mathbb{F}_{p}^{2n} \; : \; <x,y>=0, \; \forall y \in \mathbb{F}_{p}^{2n},
\; Q(x)=0\}$$ be the radical of $R(V)$. Following \cite{Ki},
define a polynomial
$$H_{n,p}(T;X):=\left\{\begin{array}{ll}
  1 & \text{if } s_p=0, \\
  \prod_{j=1}^{[(s_p-1)/2]}(1-p^{2j-1}X^2) & \text{if } s_p>0, \; s_p \text{ odd},  \\
  (1+\lambda_p(T)p^{(s_p-1)/2}X)\prod_{j=1}^{[(s_p-1)/2]}(1-p^{2j-1}X^2) &
  \text{if } s_p>0, \; s_p \text{ even},
\end{array}\right.$$
where for even $s_p$ we denote $$\lambda_p(T):=\left\{
\begin{array}{ll}
  1 & \text{if } W \text{ is a hyperbolic space or } s_p=2n, \\
  -1 & \text{otherwise}.
\end{array}\right.$$
Following \cite{Ko}, for a nonnegative integer $\mu$, define
$\rho_{T}(p^{\mu})$ by
$$\sum_{\mu \geq0}\rho_{T}(p^{\mu})X^{\mu}:=\left\{
\begin{array}{ll}
  (1-X^2)H_{n,p}(T;X), & \text{if } p|f_T, \\
  1 & \text{otherwise}.
\end{array}
\right.$$ We extend the functions $\rho_T$ multiplicatively to
natural numbers $\mathbb{N}$ by defining
$$\sum_{\mu \geq0}\rho_{T}(p^{\mu})X^{-\mu}
:=\prod_{p|f_p}((1-X^2)H_{n,p}(T;X)).
$$
Let $$\mathcal{D}(T):=GL_{2n}(\mathbb{Z})\setminus \{ G\in
M_{2n}(\mathbb{Z})\cap GL_{2n}(\mathbb{Q})\; : \; T[G^{-1}] \text{
half-integral} \},$$ where $GL_{2n}(\mathbb{Z})$ operates by
left-multiplication and $T[G^{-1}]=T'G^{-1}T$. Then
$\mathcal{D}(T)$ is finite. For $a\in \mathbb{N}$ with $a|f_T,$
 let
\begin{equation}\label{phi}
\phi(a;T):=\sqrt{a}\sum_{d^2|a}\sum_{G\in
\mathcal{D}(T),|\det(G)|=d}\rho_{T[G^{-1}]}(a/d^2).
\end{equation}
Note that $\phi(a;T) \in \mathbb{Z}$ for all $a$. With these
notations we state the following theorem:
\begin{thr}[\cite{K-K}]\label{KK}
Suppose that $g\equiv 0,1 \pmod{4}$ and let $k\in \mathbb{N}$ with
$g\equiv k \pmod{2}$. A Siegel modular form $F$ is in
$S_{k+n}^{Maass}(\Gamma_{2g})$ if and only if there exists a
modular form
$$f(z)=\sum_{n=1}^{\infty}c(n)q^n \in
S_{k+\frac{1}{2}}(\Gamma_0(4))$$ such that
$A(T)=\sum_{a|f_T}a^{k-1}\phi(a;T)c\left(\frac{|D_T|}{a^2}\right)$
for all $T$. Here, $$D_T:=(-1)^g\cdot\det(2T)$$ and
$D_T=D_{T,0}f_T^2$ with $D_{T,0}$ the corresponding fundamental
discriminant and $f_T \in \mathbb{N}$.
\end{thr}
\begin{rem}
A proof of Theorem \ref{KK} given in \cite{K-K} implies that if
$A(T)\in \mathbb{Z}$ for all $T$, then $c(m) \in \mathbb{Z}$ for all
$m\in\mathbb{N}$.
\end{rem}
\begin{proof}[Proof of Theorem \ref{main6}]
From Theorem \ref{KK} we can take
$$f(z)=\sum_{n=1}^{\infty}c(n)q^n \in
S_{k+\frac{1}{2}}(\Gamma_0(4))\cap \mathbb{Z}_{p}[[q]]$$ such that
$$F(Z)=\sum_{T>0}A(T)q^{tr(TZ)}=
\sum_{T>0}\sum_{a|f_T}a^{k-1}\phi(a;T)c\left(\frac{|D_T|}{a^2}\right)q^{tr(TZ)}.$$
By Theorem \ref{main-1}, there exists a positive integer $b$ such
that, for every positive integer $m$,
$$c(p^{b}m) \equiv 0 \pmod{p^j},$$
since $k \equiv 2 \text{ or } 3 \pmod{\frac{p-1}{2}}$. Suppose
that $p^{b+2j}||D_T|$. If $p^{j}|a$ and $a|f_T$, then
$$a^{k-1}\phi(a;T)c\left(\frac{|D_T|}{a^2}\right) \equiv 0
\pmod{p^j}.$$ If $p^j \nmid a$ and $a|f_T$, then
$p^{b}\left|\frac{|D_T|}{a^2}\right.$ and
$a^{k-1}\phi(a;T)c\left(\frac{|D_T|}{a^2}\right) \equiv 0
\pmod{p^j}.$
\end{proof}
\medskip
\section*{\bf Acknowledgement}We  thank the referee for   many
helpful comments which have improved our exposition.

\end{document}